\documentclass[onecolumn]{article}

\usepackage{amssymb}
\usepackage[english]{babel}
\usepackage{hyperref}

\newtheorem{theorem}{Theorem}[section]
\newtheorem{lemma}[theorem]{Lemma}
\newtheorem{proposition}[theorem]{Proposition}

\newtheorem{remark}{\it Remark\/}

\newcommand{\bz}{{\bf z}}

\author{Hugo Jim\'enez-P\'erez}

\begin{document}
\title{Geometrization of symplecticity conditions for implicit schemes}

\maketitle
\begin{abstract}

In this note we give simple symplecticity conditions for implicit schemes
in symplectic vector spaces. 
We consider implicit maps on generic symplectic manifold and we 
introduce the concept of \emph{consistent implicit maps} to generalize 
the symplecticity conditions on symplectic manifolds. Additionally,
we give a preliminary geometrical interpretation of those conditions.

\end{abstract}

\section{Introduction}

The most widely numerical methods for simulating Hamiltonian dynamics are 
symplectic integrators. Nowadays, there are a multitude of techniques and 
types of simplectic integrators, both explicit and implicit. 
Implicit methods, in general, are approximations which 
minimize the error but, most geometrical information is lost.
The geometry of explicit schemes
has been widely studied and it is well understood at present, which is
not the case for implicit  schemes.
The main conceptual difficulty for the implicit case, arrives from the fact that symplectic maps 
are natural operators on cotangent bundles to smooth manifolds: every 
diffeomorphism on the base manifold lifts to a symplectomorphism
on the cotangent bundle; this mapping is called the \emph{cotangent lift}  \cite{AM78,LM87}. 
Such a symplectomorphism consists on a covariant and a contravariant part acting in opposite
directions, specifically, the mapping on the fibers is defined by the inverse 
of the pullback. This fact is connected with the idea that generating functions of type {\it II}
and {\it III} are well adapted for constructing symplectic integrators, 
but not the generating functions of type {\it I} nor {\it IV}.

To be more specific, suppose that two symplectic manifolds $(M_1,\omega_1)$ and $(M_2,\omega_2)$ are given,
with canonical coordinates $(q,p)\in M_1$ and $(Q,P)\in M_2$, and a symplectomorphism
$\phi:M_1\to M_2$ between them. 
For generating functions of type {\it II} and {\it III} the interchange of geometrical information 
goes from $(q,P)\mapsto (Q,p)$ and $(Q,p)\mapsto(q,P)$ respectively. The image 
coordinates are obtained by solving the Hamilton-Jacobi equation. Looking for intermediate 
points preserving the geometrical structure, we must keep this balance, 
generalized by the rule 
\begin{eqnarray}
(\alpha q+(1-\alpha)Q, (1-\alpha)p + \alpha P)\quad \mapsto \quad ((1-\alpha)q + \alpha Q, \alpha p + (1-\alpha) P) 
  \label{eqn:alpha}
\end{eqnarray}
which is encoded in the Liouvillian form associate to the Hamiltonian system \cite{Jim15a}.
In this generalization, generating functions of type {\it II} and {\it III} correspond
to the expression (\ref{eqn:alpha}) with values $\alpha= 1$ and $\alpha= 0$ respectively.
The other well-known case, corresponding to $\alpha = 1/2$, is the 
mid-point rule. Other possibilities are not take into account since 
for $\alpha\neq 1/2$ the methods obtained are of first order and  for 
$\alpha=1/2$ is second order (the symmetric case). 

With this paper we start a systematic study on implicit symplectic 
integrators from the geometrical point of view. We outline its content. In Section 2, we state
the definitions of explicit and implicit symplectic integrators and we give 
simple simplecticity conditions for the case of linear symplectic spaces, modeled by
$(\mathbb R^{2n},\omega_0)$. They are already known results. In Section 3, 
we state some preliminary definitions in order to restrict the type of 
implicit schemes where the conditions apply. We are interested
on implicit maps wich can be 
given by the composition of two explicit maps. The main idea is to find an
intermediary point given in terms of two consecutive points of the
discretized flow, such that we can construct a flow line passing by three
successive points on the manifold. Then, we generalize and give an interpretation 
of the simplecticity conditions stated in Section 2 to the 
case of a Hamiltonian system on a generic symplectic manifold.

\section{Explicit and implicit symplectic integrators}
In a single phrase, we can state the definition of a symplectic integrator as follows: 
{\it A symplectic integrator for a Hamiltonian system is a numerical method 
which preserves the structure of the Hamiltonian vector field.}
However, the symplecticity of an integrating method only constrains the numerical scheme 
to preserve the form of the vector field. Moreover, if the Hamilton equations 
are independent of the time ({\it i.e.} if the Hamiltonian function $H$ is \emph{autonomous}), 
a symplectic integrator preserves also the energy integral.
Energy preservation restricts the numerical solution to be ``close'' to a submanifold 
$\Sigma_h= H^{-1}(h)$ of 
codimension ${\rm codim}(\Sigma_h)=1$, which is advantegeous for a lower-dimensional 
Hamiltonian system, but weak for a higher-dimensional system. Preserving the form of the vector
field says that we can apply the method to the image point without any additional analysis
but no additional constraints are given. The reason of this, rise as a consequence
of Darboux's theorem which states that the symplectic structure does not recognizes
the local structure of the Hamiltonian flow. Unfortunatelly, numerical integration 
is intrinsecally a local procedure and a new point of view is needed to go beyond 
in the subject of numerical symplectic integration.

We are conviced that implicit schemes can be very accurated methods for simulating 
the Hamiltonian flow, giving additionally the right direction of the numerical flow.
To have suitable numerical integrators we must return to the roots of the geometrical problem.
We start by some basic definitions and we state new definitions formalizing the implicit
schemes of our interest.

\subsection{Some basic definitions}
A \emph{symplectic manifold} is a $2n$-dimensional manifold $M$ equiped with 
a non-degenerated, skew-symmetric, closed 2-form $\omega$, such that
at every point $m\in M$, the tangent space 
$T_mM$, has the structure of a symplectic vector space.

A Hamiltonian system $(M,\omega,X_H)$ is a vector field
$X=X_H$ on a symplectic manifold $(M,\omega)$ such that
\begin{eqnarray}
   i_{X_H}\omega = -dH,
   \label{eqn:Ham:Def}
\end{eqnarray}
for a differentiable function $H:M\to\mathbb R$ known as the total energy or the 
Hamiltonian function. 
There is an alternative definition in vector field form as
\begin{eqnarray}
   X_H = J\nabla H,
  \qquad{\rm with\ evolution\ equations}\qquad 
    \dot \bz = J\nabla_\bz H(\bz),\qquad \bz\in M.
    \label{eqn:Field}
\end{eqnarray}
where $J$ is the canonical complex structure on $T_\bz M$ given by 
\begin{eqnarray}
    J=\left(
    \begin{array}[h]{cc}
        0_n & -I_n\\
        I_n & 0_n
    \end{array}
    \right), \qquad 0_n,I_n\in \mathbb M_{n\times n}(\mathbb R).
    \label{eqn:J}
\end{eqnarray}
and $\nabla$ is the standard gradient operator on $T_\bz M\cong \mathbb R^{2n}$.

Denoting by $\phi^t_H$ the Hamiltonian flow it is well-known that 
for each fixed $h\in\mathbb R$, the map $\phi^h_H$ is a symplectic map.
Let $\bz_0\in M$ be a point on the symplectic manifold and $\bz(t)$ the
integral curve to $X_H$ such that $\bz_0=\bz(0)$. By definition of the flow,
the mapping 
\begin{eqnarray*}
    \bz(t+h) = \phi^h_H (\bz(t))
\end{eqnarray*}
will propagates the solution from the time $t$ to time $t+h$. 
A \emph{symplectic algorithm} with stepsize $h$, is the numerical 
approximation $\psi_h$ of the 
flow $\phi^h_H:M\to M$, which is an isometry of the symplectic form
$\omega$ . 
Specifically, consider the exact solution $\bz(t)$ of a Hamiltonian system
for the time $t\in [0,T]$, and a discretization of the time interval $\{t_i\}_{i=0}^N$ such that
$t_0=0$, $t_N=T$, $h=T/N=t_{k+1}-t_k$, and $\bz_k=\bz(t_k)$ for $0\le k < N$. 

We define an \emph{explicit symplectic integrator} 
as a map 
\begin{eqnarray*}
    \psi_h: U\subset M & \to & U\\
    \bz_k & \mapsto & \bz_{k+1} = \psi_h(\bz_k)
\end{eqnarray*}
smooth on $h$ an $H$, such that $(\psi_h)^*\omega = \omega$, where 
$(\psi_h)^*$ is the pullback of $\psi_h$.

In an analogous way, we define an \emph{implicit symplectic integrator} 
as a map 
\begin{eqnarray*}
    \varphi_h: U\times U & \to & U\\
    (\bz_k, \bz_{k+1}) & \mapsto & \bz_{k+1} = \varphi_h(\bz_k, \bz_{k+1})
\end{eqnarray*}
smooth with respect to $h$ and $H$, and such that $\varphi_h^*\omega = \omega$.
However, we have a problem to state what means $\varphi_h^*\omega = \omega$ in the implicit case.
Some authors consider maps of type 
\begin{eqnarray}
	\varphi_h(\bz_k, \bz_{k+1}) = \bz_k + \phi_h(\bz_k, \bz_{k+1}).
\end{eqnarray}
or the general implicit rule  
\begin{eqnarray}
  \Psi(\bz_{k+1}, \bz_k)= \bz_{k+1} - \varphi_h(\bz_{k+1}, \bz_k)= 0
  \label{eqn:impli}
\end{eqnarray}
The test of symplecticity is obtained by implicit differentiation 
of expression (\ref{eqn:impli}), obtaining the linearized algoritm 
$\delta \bz_{k+1}=A\delta \bz_k$, where $A$ is called by some authors
the \emph{linearized amplification 
matrix} of the scheme \cite[Appx II]{SiT92}. The matrix $A$ is given by $A=A_1^{-1}A_2$ where
\begin{eqnarray*}
    A_1 = \frac{\partial \Psi(\bz_{k+1},\bz_k)}{\partial \bz_{k+1}} \quad
    {\rm and }\quad A_2 =- \frac{\partial \Psi(\bz_{k+1},\bz_k)}{\partial \bz_{k}},
\end{eqnarray*}
and the implicit mapping $\Psi$ is simplectic if $A_2^TA_1^{-T}JA^{-1}_1A_2=J$. 
With this information we can state the following result on $\mathbb R^{2n}$,
equiped with the canonical symplectic form $\omega_0$, considering $(\mathbb R^{2n},\omega_0)$
as a symplectic manifold.

\begin{proposition}\label{prop:1}
    Let $U\subset \mathbb R^{2n}$ be a convex open set and $\bz_k,\bz_{k+1}\in U$ two interior points.
    Consider a point $\bar\bz\in U$ 
    and a differentiable map $f:U\times U\to U$ such that $\bar\bz= f(\bz_k,\bz_{k+1})$.
    Denote the matrices of partial derivatives by
    \begin{eqnarray}
        B= \frac{\partial f(\bz_k,\bz_{k+1})}{\partial \bz_k} \qquad{\rm and}\qquad
        C= \frac{\partial f(\bz_k,\bz_{k+1})}{\partial \bz_{k+1}},
        \label{eqn:def:BC}
    \end{eqnarray}
    and suppose that $B$ and $C$ fulfills the following conditions
    \begin{eqnarray}
       (\imath)\ B + C = I_{2n} \qquad{\rm and}\qquad 
       (\imath\imath)\ BJ = JC^T,
        \label{eqn:mat:cond2}
    \end{eqnarray}
    where $I_{2n} \in M_{2n\times 2n}(\mathbb R)$ is the identity matrix and $J$ 
    is the almost complex structure defined in (\ref{eqn:J}). 

    Then, the map  
    \begin{eqnarray}
        \bz_{k+1} &=& \bz_k + h X_H\left( \bar\bz \right) 
        \label{eqn:sympl}
    \end{eqnarray}
    defines an implicit symplectic integrator with stepsize $h$.
    \label{teo:main}
\end{proposition}

{\it Proof.} Consider the implicit  mapping $\Psi$ given by 
\begin{eqnarray}
    \Psi(\bz_k,\bz_{k+1}) = \bz_{k+1} - \bz_k - hX_H\left(f (\bz_k, \bz_{k+1}) \right) = 0.
    \label{eqn:th1}
\end{eqnarray}
Implicit differentiation of (\ref{eqn:th1}) using  the chain rule and expressions (\ref{eqn:def:BC}) gives 
\begin{eqnarray}
        \frac{\partial \Psi(\bz_k,\bz_{k+1})}{\partial \bz_k} &=& I 
            -h JH_{zz} B \\ 
       \frac{\partial \Psi(\bz_k,\bz_{k+1})}{\partial \bz_{k+1}} &=& - I 
	    - h JH_{zz} C 
    \label{eqn:mats}
\end{eqnarray}
where $H_{zz}$ is the Hessian matrix of $H$. 
Denote the partial derivatives of $\Psi$ by 
\begin{eqnarray*}
    A_1 =  \frac{\partial \Psi(\bz_{k},\bz_{k+1})}{\partial \bz_{k}},\quad
    {\rm and }\quad 
    A_2 = - \frac{\partial \Psi(\bz_{k},\bz_{k+1})}{\partial \bz_{k+1}}.
\end{eqnarray*}
The amplification matrix of the linearized system
is $A=A_2^{-1}\circ A_1$, and $\Psi$ is symplectic if the matrix $A$ of the linearized system 
is symplectic. We recall that A is symplectic if and only if 
$A^T=A_1^{T}\circ A_2^{-T}$ is symplectic, \emph{i.e.} if the equality
$ A_2^{-1}\circ A_1 J A_1^{T}\circ A_2^{-T} = J$ holds, or equivalently if 
$ A_1JA_1^T - A_2JA_2^T = 0.$

Using the last expression, symplecticity condition becomes
\begin{eqnarray*}
    (I-hJ H_{zz}B )J(I-hJH_{zz} B)^T - (I+hJH_{zz} C)J(I+hJ H_{zz} C)^T=0.
\end{eqnarray*}
Developping and symplifying we have
\begin{eqnarray*}
    h\left( JH_{zz} B J + JB^TH_{zz}^TJ^T + JH_{zz} CJ +JC^TH_{zz}^TJ^T \right)\qquad\qquad\\
     - h^2 \left(JH_{zz} BJB^TH_{zz}^TJ^T -JH_{zz} CJC^TH_{zz}^TJ^T \right)=0.
\end{eqnarray*}
Using the facts that $H_{zz}=H_{zz}^T$, $J^T=-J$ and $h\neq0$, factoring and reordering 
we obtain the system of equations 
\begin{eqnarray*}
   0&=& H_{zz} ( B +  C) - (B^T + C ^T)H_{zz}\\
   0 &=&  BJB^T- CJC^T .
\end{eqnarray*}
First equation is satisfied applying hypothesis $(\imath)$.
For the second equation we substitute $B^T=I-C^T$ in $BJB^T-CJC^T$ to obtain succesively
\begin{eqnarray*}
    BJB^T-CJC^T= BJ(I-C^T)-CJC^T = BJ - (B+C)JC^T \stackrel{(\imath)}{=} BJ - JC^T.
\end{eqnarray*}
By hypothesis $(\imath\imath)$ the second equation is satisfied. Consequently, the implicit method
(\ref{eqn:sympl}) is symplectic as we want to prove.
$\hfill\square$

\begin{lemma}
  \label{lem:equiv:ham}
  Let  $B,C\in GL(2n)$ be two matrices defined on a symplectic vector space $(V,\omega_0)$ and 
  $J\in GL(2n)$ the complex matrix associated to $\omega_0$. If $B+C= I$ is the identity matrix,
  then the following statements are equivalents:
  \begin{enumerate}
      \item $B J = JC^{T}$,
      \label{eqn:impl:ham}
      \item $(B - C)$ is a Hamiltonian matrix.
  \end{enumerate}
\end{lemma}
{\it Proof.-} We recall that a square matrix $A\in GL(2n)$ is Hamiltonian if $A^TJ + JA= 0$.
A direct computation shows that 
\begin{eqnarray*}
   (B-C)J+J(B-C)^T &=& BJ - CJ + JB^T - JC^T\\
       &=& BJ - (I-B)J +J(I-C^T) - JC^T\\
       &=& 2(BJ - JC^T)\\
       &=& 0
\end{eqnarray*}
where we used the fact that $A^T$ is Hamiltonian if and only if $A$ does.

$\hfill\square$

Using this lemma, we have that the map (\ref{eqn:sympl}) is symplectic when 
the matrices $B$ and $C$ satisfy that their addition is the 
identity matrix and their difference is a Hamiltonian one. 
From the previous results, they can be rewritten as 
\begin{eqnarray}
   B:= \frac{1}{2}(I+b),\qquad C:= \frac{1}{2}(I-b),
\end{eqnarray}
and conditions $(\imath)$ and $(\imath\imath)$ become  $b^TJ+Jb=0$.
\begin{remark}
In a slightly different context, Ge and Dau-liu obtain a similar condition for some matrix $b$
(presumably due to Feng) when looking for examples of 
generating functions which are invariant under symplectic transformations of Feng's type 
$\alpha_0$ (see \cite[sec. 6]{GD95}). 
In particular, the case $b=0$ corresponds to the symplectic midpoint rule which 
Feng associated to the Poincar\'e's generating function and the study of its invariance under symplectic 
transformations is due to Weinstein  \cite{Wei72}. However, in \cite{Jim15f} the author shows
that Poincar\'e's generating function does not produce a symplectic map profitable for 
numerical integrators. It looks that this condition is related to a different type of
symplectic maps adapted for dealing with periodic orbits.
\end{remark}

Note that the matrices $B$ and $C$ are well defined by the natural diffeomorphisms
\begin{eqnarray}
   T^*\mathbb R^n\cong T\mathbb R^n \cong \left(\mathbb R^{2n}\right)^* \cong \mathbb R^{2n}.
\end{eqnarray}

However, for a generic symplectic manifold $(M,\omega)$, $B$ and $C$ will be linear operators 
acting on different linear spaces (for instance, they act on different fibers of the tangent bundle). 
In the next section we develop the equivalent conditions for the symplectic generic case.

\section{Consistent implicit maps}
We are looking for a geometrical generalization to the conditions of Proposition 
(\ref{prop:1}) when the phase space is considered as a generic smooth manifold
of dimension $2n$.
In this case, tangent vectors on generic curves belong to different tangent spaces 
and we study what is the generalization of the matrices $B$ and $C$ as objects of the differential 
geometry of the smooth manifold. In what follows, we consider $M$ as any smooth manifold of 
arbitrary dimension and 
$U\subset M$ an open convex set of $M$. Restrictions on the dimension and geometry of $M$
will be stated when necessary. All the results studied here are localized into the open
set $U\subset M$.

Let $\phi:U\times U\to U$ an implicit map such that $\bz_{k+1}= \phi(\bz_k,\bz_{k+1})$. 
We say that $\phi$ is \emph{consistent} if there exists $\bar \bz\in U$ and two 
local diffeomorphisms $\psi_1,\psi_2:U\to U$ with $\bar \bz = \psi_1(\bz_k)$ 
and $\bar \bz = \psi_2(\bz_{k+1})$, such that:
\begin{enumerate}
 \item it is possible to rewrite $\phi$ in the form 
   \begin{eqnarray}
      \phi(\bz_k,\bz_{k+1})=  \psi_2^{-1}(\bar \bz)=\bz_{k+1},
   \end{eqnarray}
  \item the limit 
   \begin{eqnarray}
      \lim_{\bz_{k+1}\to\bz_k} \psi_i = Id, \qquad i=1,2.
   \end{eqnarray}
   hold.
\end{enumerate}
We call $\bar\bz$ the \emph{consistency point}.

There is a natural 
local diffeomorphism $\psi:U\to U$ given by $\psi=\psi_2^{-1}\circ\psi_1$ 
which is the explicit counterpart of the implicit map $\phi$. This is 
called the \emph{consistency map} and it enables the construction of solutions 
passing by the three points $\bz_k$, $\bar \bz$ and $\bz_{k+1}$. 

\begin{lemma}
  \label{lem:mapp:exists}
  For every consistent implicit map $\phi:U\times U\to U$ it is possible to generate an implicit map 
  $\rho: U\times U\to U$ such that 
  \begin{eqnarray}
    \bar \bz = \rho(\bz_{k}, \bz_{k+1}).
  \end{eqnarray}
\end{lemma}
{\it Proof.} Since $\phi$ is a consistent implicit map, there exist local diffeomorphisms 
$\psi_1,\psi_2\in {\rm Diff}_0(U)$ and $\bar\bz\in U$ such that $\bar \bz= \psi_1(\bz_n)$ and 
$\bar \bz= \psi_1(\bz_{n+1})$.
Consider a convex combination 
\begin{eqnarray}
   \rho(\bz_k, \bz_{k+1}): = a\psi_1(\bz_k)+(1-a)\psi_2(\bz_{k+1}), \qquad a\in \mathbb R.
   \label{eqn:map}
\end{eqnarray}
Then we have successively 
\begin{eqnarray*}
   \rho(\bz_k, \bz_{k+1}) &=& a\psi_1(\bz_k)+(1-a)\psi_2(\bz_{k+1}),\\
      &=& a\bar\bz + (1-a)\bar\bz\\
      &=& \bar\bz.
\end{eqnarray*}
$\hfill\square$

Note that the general case $a\in \mathbb R$ is well defined, however, we are looking
for localized maps in the open $U$. Moreover, we want to constrain the point $\bar\bz$
to be an intermediate point on the same flow line of $\bz_k$ and $\bz_{k+1}$ then we 
restrict its domain to $a\in[0,1]$. From now on, we will see the map (\ref{eqn:map})
as a partition of the unity.

\begin{lemma}
  \label{lem:mapp:id}
  If $\phi:U\times U\to U$ is a consistent implicit map and $\rho:U\times U \to U$ the map
  given by (\ref{eqn:map}) with image on the consistency point $\bar \bz = \rho(\bz_{k}, \bz_{k+1})$.
  Then its tangent
  map $T\rho : TU\times TU\to TU$ corresponds to the identity map on $T_{\bar\bz}U$.
\end{lemma}
{\it Proof.} By the consistency hypothesis, we have a point $\bar\bz\in U$ and two local
diffeomorphisms $\psi_1,\psi_2\in {\rm Diff}_0(U)$ satisfaying the hypothesis of Lemma 
\ref{lem:mapp:exists}. For every vector $v\in T_{\bar\bz}U$, 
there exist vectors $v_1\in T_{\bz_k}U$ and $v_2\in T_{\bz_{k+1}}U$ such that 
$v = T\psi_1(v_1)$ and $v = T\psi_2(v_2)$, where $T\psi_i:TU\to TU$ are the 
tangent maps to $\psi_i$ for $i=1,2$.

To be more specific, we have
\begin{eqnarray}
   T\psi_1|_{\bz_k}: T_{\bz_k}U\to T_{\bar\bz}U\quad{\rm and}\quad T\psi_2|_{\bz_{k+1}} : T_{\bz_{k+1}}U\to T_{\bar\bz}U.
\end{eqnarray}
If $\rho:U\times U \to U$ is of the form (\ref{eqn:map}), its tangent map 
$T\rho : T_{\bz_k}U\times T_{\bz_{k+1}}U\to T_{\bar\bz}U$ takes 
$(v_1,v_2)^T\mapsto v=a T\psi_1(v_1) + (1-a)T\psi_2(v_2)$.

Finally note that 
\begin{eqnarray}
   T\rho \left(
      \begin{array}{c}
         v_1\\
         v_2
      \end{array}
      \right)
     = 
     \left( aT\psi_1, (1-a)T\psi_2 \right)
     \left(
      \begin{array}{c}
         v_1\\
         v_2
      \end{array}
      \right)
      = v
\end{eqnarray}
where we used the identity map $v=Id(v)$ on $T_{\bar\bz}U$ on the right hand side.

$\hfill\square$

In local coordinates we obtain the expression
    \begin{eqnarray}
      \frac{\partial \rho(\bz_k,\bz_{k+1})}{\partial \bz_k}(v_1) + 
      \frac{\partial \rho(\bz_k,\bz_{k+1})}{\partial \bz_{k+1}}(v_2) = I_{2n}(v),
      \label{eqn:consistent}
  \end{eqnarray}
which corresponds to the generalization of $B+C = I_{2n}$.

\subsection{Symplectic constraints for consistent implicit maps}
From now on, the analysis concerns the symplectic case and the manifold of interest are the generic
symplectic manifold $(M,\omega)$ of dimension ${\dim( M)}=2n$.

We say that a consistent implicit map $\phi:U\times U\to U$ in an open convex set of a symplectic manifold 
$(M,\omega)$ \emph{interleaves} a symplectic map if its consistency map $\psi$ is 
symplectic. 

\begin{lemma}
  \label{lem:mapp:symp}
   Let $\phi:U\times U\to U$ be a consistent implicit map and $\psi=\psi_2^{-1}\circ\psi_1$ its 
   consistency map. 
  Then $\phi:U\times U\to U$ interleaves a symplectic
  map if 
  \begin{eqnarray}
      (T\psi_1)^{-T}J(T\psi_1)^{-1} = (T\psi_2)^{-T}J(T\psi_2)^{-1}.
      \label{eqn:impl:symp}
  \end{eqnarray}
  where $T\psi_1$ and $T\psi_2$ are the tangent maps and $(\cdot)^{-T}=\left(\left(\cdot\right)^{-1}\right)^T$ 
  is the transpose of the inverse map. Moreover, if $\psi_1$ and $\psi_2$ are symplectic 
  then 
  \begin{eqnarray}
      (T\psi_1)^TJT\psi_1 = (T\psi_2)^TJT\psi_2.
      \label{eqn:impl:symp2}
  \end{eqnarray}
\end{lemma}
{\it Proof.} By definition, $\phi$ interleaves a symplectic map if its 
consistency map $\psi=\psi_2^{-1}\circ\psi_1$ is symplectic, it means if 
$\psi^*\omega = \omega$. This holds if and only if we have succesively the following equalities
\begin{eqnarray*}
   \psi^*\omega &=& \left(\psi_2^{-1}\circ \psi_1\right)^*\omega\\
         &=& \psi_1^*\circ\left(\left(\psi_2^{-1}\right)^*\omega\right)\\
         &=& \omega
\end{eqnarray*}
from where we obtain $\left(\psi_2^{-1}\right)^*\omega  = \left(\psi_1^{-1}\right)^* \omega$.
This condition is equivalent to the equation 
  \begin{eqnarray}
    \left( T \psi_1^{-1} \right)^T  J 
    \left( T \psi_1^{-1}\right) =
    \left( T \psi_2^{-1}\right)^T  J  
    \left( T \psi_2^{-1}\right).
  \end{eqnarray}
on $T_{\bar\bz}U$ for $\bar\bz\in U$ the consistency point of $\phi$. This proves the first 
result. For the second part, the symplectic hypothesis on $\psi_1$ and $\psi_2$ implies
that $\psi_1^{-1}$ and $\psi_2^{-1}$ are also symplectic (local) diffeomorphisms by the group 
property of $Sp(U,\omega)$. 

$\hfill\square$

\begin{remark}

Condition (\ref{eqn:impl:symp}) says nothing about the symplecticity of the mappings $\psi_1$ and $\psi_2$, but
only that the composition is symplectic. 
For instance, suppose that 
\begin{eqnarray}
  \psi_1 = f\circ g \quad {\it and} \quad \psi_2 = f\circ h, 
\end{eqnarray}
for $h,g\in Sp(U,\omega)$, $h\neq g$ and 
$f\in {\rm Diff}_0(U)$, but $f\notin  Sp(U,\omega)$. Then, $\psi_1,\psi_2\notin Sp(U,\omega)$ however 
$\psi_2^{-1}\circ \psi_1 = h\circ g\in Sp(U,\omega)$.
\end{remark}

In other words, condition (\ref{eqn:impl:symp}) says that the complex structures 
$J$ on $T_{\bz_k}U$ and $T_{\bz_{k+1}}U$ 
are equivalent, but it might not be well defined on $T_{\bar\bz}U$ when $\psi_1,\psi_2\notin Sp(U,\omega)$. 
Expression (\ref{eqn:impl:symp2}), on the other hand, says that there exists a structure 
on $T_{\bar\bz}U$ which is equivalent to those $J$'s on $T_{\bz_k}U$ and $T_{\bz_{k+1}}U$ and
it defines a complex structure on $T_{\bar\bz}U$.

Since we are interested in simulating Hamiltonian flows, we impose the condition that 
the components of the consistency map $\psi_1,\psi_2\in Sp(U,\omega)$ are symplectic in the 
rest of this work.

\begin{lemma}
  \label{lem:mapp:ham}
  Let  $\phi:U\times U\to U$ be a consistent implicit map interleaving a symplectic map.
  If the components of the consistency map $\psi=\psi_2^{-1}\circ \psi_1$ are symplectic, then 
  \begin{eqnarray}
      T\psi_1 - T\psi_2 : TU\times TU\to TU
  \end{eqnarray}
   is a Hamiltonian operator on $T_{\bar\bz}U$.
\end{lemma}
{\it Proof.-} {\it 1)}
Consider the map $\bar\bz = \rho(\bz_{n},\bz_{n+1}) = a \psi_1(\bz_n) + (1-a)\psi_2(\bz_{n+1})$
as a partition of the unity of two local charts associated to $\psi_1$ and $\psi_2$.
Since $\psi_1,\psi_2\in Sp(U,\omega)$ we can consider the curve 
\begin{eqnarray}
   \gamma_\tau = \tau\psi_1 + (1 - \tau)\psi_2
   \label{eqn:curve0}
\end{eqnarray}
as a one parameter family of symplectic diffeomorphisms with image $\bar\bz$.
Now, we consider the properties of $T_{\bar\bz} U$ as a vector space.  

For every $v\in T_{\bar\bz}U$,
and for every $\tau\in[0,1]$, there exist $\hat{v}_1,\hat{v}_2\in T_{\bar\bz}U$ such that
$v=\tau \hat{v}_1 + (1-\tau)\hat{v}_2$. This defines a $(2n-1)$-dimensional subspace of $T_{\bar\bz}U$
(the hyper-plane perpendicular to $v$). We consider the space of smooth curves joining the origin 
in $T_{\bar\bz}U$ with the end point of $v$, and it is possible to write 
\begin{eqnarray}
   \gamma_\tau = \tau\hat{v}_1 + (1- \tau)\hat{v}_2
   \label{eqn:curve1}
\end{eqnarray}
We look for vectors $v_1\in T_{\bz_n}U$ and $v_2\in T_{\bz_{n+1}}U$
such that $\hat{v}_1=T\psi_1(v_1)$ and $\hat{v}_2=T\psi_2(v_2)$, and the curve (\ref{eqn:curve1}) can be
rewritten as
\begin{eqnarray}
   \gamma_\tau = \tau T\psi_1(v_1) + (1 -\tau)T\psi_2(v_2) 
   \label{eqn:curve2}
\end{eqnarray}
Since $\psi_1,\psi_2\in Sp(U,\omega)$ we can consider the curve (\ref{eqn:curve0}) 
as a one parameter family of symplectic diffeomorphisms with image $\bar\bz$,
joining the symplectomorphisms $\psi_1$ and $\psi_2$.

Taking the derivative with respect to the parameter $\tau$ we have
\begin{eqnarray}
   \frac{ \partial \gamma_\tau}{\partial \tau} = T\psi_1 -T\psi_2
\end{eqnarray}

Which by definition is a Hamiltonian operator.
$\hfill\square$

\begin{remark}
Other parameterizations like $\bar\gamma_\tau=\cos^2(\tau)\hat{v}_1 + \sin^2(\tau)\hat{v}_2$
with derivative
\begin{eqnarray}
   \frac{ \partial \bar\gamma_\tau}{\partial \tau} = 2\cos(\tau)\sin(\tau)\left(T\psi_1 -T\psi_2\right)
\end{eqnarray}
or $\bar{\bar\gamma}_\tau= {\rm cn}^2(\tau,k)\hat{v}_1 + {\rm sn}^2(\tau, k)\hat{v}_2$, $k\in(0,1) $
with derivative
\begin{eqnarray}
   \frac{ \partial \bar{\bar\gamma}_\tau}{\partial \tau} =2{\rm sn}(\tau,k){\rm cn}(\tau,k){\rm dn}(\tau,k)
	\left( T\psi_1 -T\psi_2\right)
    \label{eqn:par2}
\end{eqnarray}
leads to equivalent results modulo a scalar function on the 
open set $(0,1)$. All of them are parallel 
Hamiltonian operators.
Expression (\ref{eqn:par2}) is given in terms of the elliptic functions of Jacobi
${\rm sn}(\tau,k)$, ${\rm cn}(\tau,k)$ and ${\rm dn}(\tau,k)$, which also satisfy 
${\rm cn}^2(\tau,k) + {\rm sn}^2(\tau,k) = 1$
\end{remark}

The conditions $\imath) B+C=I$, and $\imath\imath)BJ=JC^T$ 
in Proposition \ref{prop:1} only constraint the implicit mapping to be consistent and
symplectic, but no relationship with any particular Hamiltonian flow is stated. 
From the geometrical interpretation of condition $\imath\imath)$ the only requirement
was that vectors $\hat{v}_1$ an $\hat{v}_2$ must 
depend smoothly on the parameter $\tau$ satisfied by the use of smooth curves. However, we are
interested in constrain these vectors  such that they 
were related with the Hamiltonian vector field in the following way: 
$\hat{v}_1= T\psi_1(X_H(\bz_{n}))$ and $\hat{v}_2= T\psi_2(X_H(\bz_{n+1}))$.
A good candidate might be $$v=\frac12 \left(T\psi_1(X_H(\bz_n))+ T\psi_2(X_H(\bz_{n+1}))\right),$$
however, since we know the form of the Hamiltonian vector field, we have that 
$v = X_H(\bar\bz)$. The problem is not to find the good vector $v$, instead 
we look for a point $\bar\bz\in U$ on the Hamiltonian flow producing the good 
value for $v=X_H(\bar\bz)$.

\begin{remark}
   Conditions stated in Proposition \ref{prop:1} give rise to continuous high-dimensional
   families of implicit symplectic integrators. We study some families in \cite{Jim15a}
   using the classical framework used in the method of generating functions.
\end{remark}

%
%


The proof of Lemma \ref{lem:mapp:ham} looks a slightly tricky and artificial, at the same 
time the generalization of the relation $B+C=I_{2n}$ might be unsatisfactory. This
arises since the natural definitions for symplectic and Hamiltonian operators are given 
between manifolds of the same dimension.
   
This problem is solved if we define all the operators on the product manifold  
of two copies of the  symplectic manifold $(M,\omega)$. For instance, consider the manifold 
$\tilde M = M_1\times M_2$ with the 2-form
$\omega_\ominus = \pi^*_1\omega_1 - \pi^*_2\omega_2$, defined by the pullback of 
the canonical projections $\pi_i:\tilde M\to M_i$, $i=1,2$; the couple 
$(\tilde M,\omega_\ominus)$ is a symplectic 
manifold of dimension $4n$ (see \cite{AM78}). 
The procedure consists in create a symplectic path in $\tilde M$ 
from $\left( M_1\times \{0\}\right)$ to $\left(\{0\} \times M_2\right)$.
It means, a continuous family of symplectic subspaces in $\tilde M$ of dimension $2n$
joining $M_1$ to $M_2$. For some intermediary element in this family, we project 
the (mixed) coordinates on the original manifold $(M,\omega)$, and these becomes
the coordinates of the point $\bar \bz\in U\subset M$ that we are looking for.
However, there are some subtleties which are worked out in the companion article \cite{Jim15a}.

\section*{Acknowledgements}
The author thanks J.P. Vilotte and B. Romanowicz for their support and constructive
criticism on this work. 
This research was developed with support from the \emph{Fondation du Coll\`ege de 
France} and \emph{Total} under the research convention PU14150472, as well as the ERC Advanced Grant 
WAVETOMO, RCN 99285, Subpanel PE10 in the F7 framework.

\end{document}